\documentclass{article}
\usepackage{CJK,CJKnumb,CJKulem,times,dsfont,ifthen,mathrsfs,latexsym,amsfonts, color}
\usepackage{amsmath,amsthm,makeidx,fontenc,amssymb,bm,graphicx,psfrag,listings, curves,extarrows}
\usepackage[top=1.0in, bottom=1.0in, left=1.25in, right=1.25in]{geometry}
\makeindex
\newtheorem{definition}{Definition}[section]
\newtheorem{theorem}{Theorem}[section]
\newtheorem{lemma}{Lemma}[section]
\newtheorem{corollary}{Corollary}[section]
\newtheorem{proposition}{Proposition}[section]
\newtheorem{remark}{Remark}[section]

\newcommand{\bt}{\begin{theorem}}
	\newcommand{\et}{\end{theorem}}
\newcommand{\bl}{\begin{lemma}}
	\newcommand{\el}{\end{lemma}}
\newcommand{\bd}{\begin{definition}}
	\newcommand{\ed}{\end{definition}}
\newcommand{\bc}{\begin{corollary}}
	\newcommand{\ec}{\end{corollary}}
\newcommand{\bp}{\begin{proof}}
	\newcommand{\ep}{\end{proof}}
\newcommand{\bx}{\begin{example}}
	\newcommand{\ex}{\end{example}}
\newcommand{\bi}{\begin{exercise}}
	\newcommand{\ei}{\end{exercise}}
\newcommand{\bo}{\begin{proposition}}
	\newcommand{\eo}{\end{proposition}}
\newcommand{\br}{\begin{remark}}
	\newcommand{\er}{\end{remark}}
\newcommand{\be}{\begin{equation}}
\newcommand{\ee}{\end{equation}}
\newcommand{\ba}{\begin{align}}
\newcommand{\ea}{\end{align}}
\newcommand{\bn}{\begin{enumerate}}
	\newcommand{\en}{\end{enumerate}}
\newcommand{\bg}{\begin{align*}}
\newcommand{\bcs}{\begin{cases}}
\newcommand{\ecs}{\end{cases}}

\newcommand{\bean}{\begin{eqnarray*}}
\newcommand{\eean}{\end{eqnarray*}}

\numberwithin{equation}{section}

\begin{document} 
\begin{CJK*}{GBK}{song}

\title{\bf {Existence of solutions for a class of nonlinear Choquard equations with critical growth}
	\thanks{E-mails: aoy15@mails.tsinghua.edu.cn}}

\date{}
\author{{\bf  Yong Ao}\\
	\footnotesize {\it Department of Mathematical Sciences, Tsinghua University,  Beijing 100084, China}}

\maketitle

\vskip0.7in

\begin{center}
\begin{minipage}{120mm}

\begin{center}{\bf Abstract}\end{center}
In this paper, we consider the nonlinear Choquard  equations of the form 
$$-\Delta u+ u=(I_{\alpha}*|u|^p)|u|^{p-2}u+|u|^{q-2}u~~~\mbox{in}~\mathbb{R}^N,$$ 
where $I_{\alpha}$ is a Riesz potential,~$\alpha\in(0,N)$,~$N\geqslant 4$,~$ p= \dfrac{N+\alpha}{N-2}$, ~$2< q < 2^{*}$. We show the existence of  solutions of the equations.

\vskip0.1in


\end{minipage}
\end{center}
\vskip0.27in

\section{Introduction}
In this paper we consider the Choquard problem 
\begin{equation}
-\Delta u+ u=(I_{\alpha}*|u|^p)|u|^{p-2}u+|u|^{q-2}u~~~\mbox{in}~\mathbb{R}^N,
\end{equation}
where $I_{\alpha}$ is the Riesz potential defined by 
$$I_{\alpha}(x)=\frac{\varGamma(\frac{N-\alpha}{2}) }{\varGamma(\frac{\alpha}{2})\pi^{N/2}2^{\alpha}|x|^{N-\alpha} }\stackrel{\triangle}{=}\dfrac{\bar{C}}{|x|^{N-\alpha}},$$
$$(I_{\alpha}*|u|^p)(x)=\int_{\mathbb{R}^N}\dfrac{\bar{C}|u(y)|^p}{|x-y|^{N-\alpha}}dy,$$
and $\varGamma $ is the Gamma function.

As is known that the solutions can be obtained by finding the critical points of the following functional 
$$I(u)=\dfrac{1}{2}\int_{\mathbb{R}^N}|\nabla u|^{2}dx+\dfrac{1}{2}\int_{\mathbb{R}^N}| u|^{2}dx-\frac{1}{2p}\int_{\mathbb{R}^N} (I_{\alpha}*|u|^p)|u|^pdx-\dfrac{1}{q}\int_{\mathbb{R}^N}|u|^qdx.$$
The following problem 
\begin{equation}
-\Delta u+ u=(I_{\alpha}*|u|^p)|u|^{p-2}u~~~\mbox{in}~\mathbb{R}^N,
\end{equation}
where  $0 <\alpha <N$, ~$ \dfrac{N+\alpha}{N} \leqslant p \leqslant\dfrac{N+\alpha}{N-2}$, has been studied in the past few years.
Problem (1.2) is a nonlocal one due to the existence of the nonlocal nonlinearity. It arises in various fields of mathematical physics, such as quantum mechanics, physics of laser beams, and the physics of multiple-particle systems. When $N=3,~~\alpha=p=2$, (1.2) turns to be the well-known Choquard-Pekar equation:
\begin{equation}
-\Delta u+u=(I_2*|u|^2)u~~~\mbox{in}~\mathbb{R}^3,
\end{equation}
which was proposed as early as in 1954 by Pekar to describe the quantum machanics of a polaron at rest in [25], and by a work of Choquard in 1976 in a certain approximation to Hartree-Fock theory for one-component plasma, see [14]. (1.3) is also known as the nonlinear stationary Hartree equation since if $u$ solves (1.3) , then $\psi (t,x)=e^{it}u(x)$ is a solitary wave of the following time-dependent Hartree equation
$$i\psi_t=-\Delta \psi- (I_2*|\psi|^p)|\psi|^{p-2}\psi,$$
see [19].
There have been many papers considering problems (1.2) and (1.3) by variational methods and also by ordinary differential equations techniques, see [8,14,16,18-21]. In [21], Moroz and Schaftingen proved that (1.2) has a nontrivial solution if and only if $ \dfrac{N+\alpha}{N} < p <\dfrac{N+\alpha}{N-2}$ (see also [11,17]), and the upper critical exponent $\dfrac{N+\alpha}{N-2}$ appears as a natural extension of the critical Sobolev exponent $\dfrac{2N}{N-2}$. 

Physical models in which particles are under the influence of an external electric field, lead to study Choquard equations in the form 
\begin{equation}
-\Delta u+ V(x)u=(I_{\alpha}*F(u))f(u)~~~\mbox{in}~\mathbb{R}^N,
\end{equation}
where $F(s)$ is the primitive of $f(s)$. When $f(s)=|s|^{p-2}s$, the authors in [4] have proved the existence of multi-bump solutions for (1.4) with deepening potential well $V(x)=\lambda a(x) + 1$ in $\mathbb{R}^3$. The existence of a nontrivial positive solution for the lower critical exponent $p=\dfrac{N+\alpha}{N}$ has been studied in [22]. And for the upper critical exponent $p=\dfrac{N+\alpha}{N-2}$, when $V(x)=\lambda$, the authors in [10] have proved the existence of a nontrivial solution in a bounded domain. For a more general nonlinearity, when $V(x)=1$, Moroz and Schaftingen have proved the exitence of the groundstate solution in [23], and the authors in [9] have proved the existence of groundstates for nonlinear fractional Choquard equations. In [13], Li and Ye have proved the existence of positive solutions with prescribed $L^2$-norm for problem (1.4). The authors in [3] have proved the existence of a positive solution when $V(x)$ is a radial function and vanishes at infinity for $p=\dfrac{N+\alpha}{N-2}$. When $\Delta$ is replaced by $\varepsilon^2\Delta$, the existence and concentration of groundstate solutions have been proved in [2].

For problem (1.1), the existence of solutions has been proved when $N=3,~0<\alpha<1,~p=2$ and $4\leqslant q<6$ in [6]. 
When $N=2$, the Riesz potential $I_\alpha$ is replaced by two-dimensional Newtonian potential $W$ defined for $x \in \mathbb{R}^2\setminus \{0\}$ by $W(x)=-\log(|x|)/(2\pi)$, then the problem has been studied for $p=2$ and $q>2$ in [7].

In our paper, we consider the  case of $N\geqslant 4$,~$ p= \dfrac{N+\alpha}{N-2}$, ~$2< q < 2^{*}$. 
 
\bt\label{th}Let $0 <\alpha <N$,~$ p= \dfrac{N+\alpha}{N-2}$, when ~$N\geqslant 5$,~$2< q < 2^{*}$, or $N=4$, $3<q<4$, then (1.1) has a nontrival solution.
\et

Consider the following minimizing problem 
\begin{equation}
A=\inf\{\frac{1}{2}\int_{\mathbb{R}^N}|\triangledown u|^{2}dx:u\in H^1(\mathbb{R}^N),~\mbox{and}~ \int_{\mathbb{R}^N}G(u)dx=1\},
\end{equation}
where $$\int_{\mathbb{R}^N}G(u)dx=\frac{1}{2p}\int_{\mathbb{R}^N} (I_{\alpha}*|u|^p)|u|^pdx+\dfrac{1}{q}\int_{\mathbb{R}^N}|u|^qdx-\dfrac{1}{2}\int_{\mathbb{R}^N}| u|^{2}dx.$$
In [5], the authors studied the following nonlinear scalar field equation 
\begin{equation}
-\Delta u=g(u)~~~in~\mathbb{R}^N,
\end{equation}
where $g: \mathbb{R}\longrightarrow \mathbb{R}$ is  a continuous function. They obtained the existence of groundstate solution of the problem by considering the above minimization problem
where $G:\mathbb{R}\longrightarrow \mathbb{R}$ is the primitive of a function $g: \mathbb{R}\longrightarrow \mathbb{R}$ which has a subcritical growth. In [12], the authors showed that the mountain pass level of the energy functional associated to (1.6) is a critical value and corresponds to the groundstate found in [5]. In [1,28], the authors studied the existence  of the least energy solutions for the critical growth. 
\vspace{0.4cm}
\bt\label{th}Let $0 <\alpha <N$,~$ p= \dfrac{N+\alpha}{N-2}$, when ~$N\geqslant 5$,~$2< q < 2^{*}$, or $N=4$, $3<q<4$, then (1.5) has a minimizer.
\et
\vspace{0.4cm}
\br\label{re} In our paper we can't obtain the least energy solution by scaling the minimizer due to the nonlocal term and the subcritical term. \er

The paper is organized as follows. In section 2,  we give some notations and prove theorem 1.1. In section 3, we conclude the proof of theorem 1.2.

\vspace{0.4cm}
\section{The Proof of Theorem 1.1}

First we give some notations. 
Denote $$A(u)= \int_{\mathbb{R}^N}|\nabla u|^{2}dx,~~~B(u)=\int_{\mathbb{R}^N} (I_{\alpha}*|u|^p)|u|^pdx,$$
$$C(u)=\int_{\mathbb{R}^N}|u|^qdx ,~~~D(u)=\int_{\mathbb{R}^N}| u|^{2}dx.$$
We consider problem (1.1) on the space$H_{rad}^1(\mathbb{R}^N)$ and denote the Nehari manifold 
$$\mathcal{N}=\{u\in H_{rad}^1(\mathbb{R}^N)\setminus \{0\}: J(u)\stackrel{\triangle}{=}\left< I'(u),u \right> =A(u)+D(u)-B(u)-C(u)=0\}, $$
and we know that the critical point of the functional must lie in the Nehari manifold. To obtain the existence of a solution of the equation, we consider the following constrained minimizing problem 
$$c=\inf\limits_{u\in \mathcal{N}}I(u).$$ 
In what follows, $S>0$ denotes the best constant of Sobolev embedding
$$D^{1,2}(\mathbb{R}^N)\hookrightarrow L^{2^*}(\mathbb{R}^N),$$
and $C$ denotes positive constant which may be different in different places.
It is important to recall the  Hardy-Littlewood-Sobolev inequality which will be frequently used in our paper, see [15].
\vspace{0.4cm}
\bo\label{p} Let $s,r>0$ and $0<\lambda <N$ with $1/s+\lambda /N+1/r=2$. If $f\in L^s(\mathbb{R}^N)$ and $g\in L^r(\mathbb{R}^N)$, then there exists a sharp constant $C(s,N,\lambda)$ such that 
$$\int_{\mathbb{R}^N}\int_{\mathbb{R}^N}\dfrac{f(x)g(y)}{|x-y|^\lambda}dxdy\leqslant C(s,N,\lambda)\lVert f\rVert_s\lVert g\rVert_r.$$
Moreover, if $s=r=2N/(2N-\lambda)$, then 
$$C(s,N,\lambda)=C(N,\lambda)=\pi^{\frac{\lambda}{2}}\dfrac{\Gamma(\frac{N-\lambda}{2})}{\Gamma(\frac{2N-\lambda}{2})}\big(\dfrac{\Gamma(\frac{N}{2})}{\Gamma(N)}\big)^{-\frac{N-\lambda}{N}},$$
and the sharp constant is achieved if and only if $g\equiv (const.)f$ and 
$$f(x)=A(\gamma^2+|x-a|^2)^{-(2N-\lambda)/2},$$
for some $A\in \mathbb{C}$, $0\neq \gamma \in \mathbb{R}$ and $a\in \mathbb{R}^N$.\eo 
Then we have 
$$\int_{\mathbb{R}^N} (I_{\alpha}*|u|^p)|u|^pdx\leqslant C_0 \lVert u \rVert_{2^*}^{2p},$$
where $$C_0=C(N,\lambda)\bar{C}=(4\pi)^{-\frac{\alpha}{2}}\dfrac{\Gamma(\frac{N-\alpha}{2})}{\Gamma(\frac{N+\alpha}{2})}\big(\dfrac{\Gamma(\frac{N}{2})}{\Gamma(N)}\big)^{-\frac{\alpha}{N}}.$$

\vspace{0.4cm}
\bl\label{lemma} The Nehari manifold $\mathcal{N}$ is not empty.\el
\vspace{0.4cm}
\noindent {\bf Proof. } For $\forall u\in H_{rad}^1(\mathbb{R}^N)\setminus \{0\}$, and $t>0$, 
$$J(tu)=t^2[A(u)+D(u)]-t^{2p}B(u)-t^qC(u),$$
then 
\begin{equation}
\begin{split}
\dfrac{d}{dt}J(tu)&=2t[A(u)+D(u)]-2pt^{2p-1}B(u)-qt^{q-1}C(u)\\
&=t[2A(u)+2D(u)-2pt^{2p-2}B(u)-qt^{q-2}C(u)]\\
&\stackrel{\triangle}{=}tg(t).
\end{split}
\end{equation}
It is easy to verify that there exist a unique $t_1>0$ such that 
$$
g(t) \begin{cases} >0 & for~~ 0<t<t_1,\\
<0  & for~~ t>t_1.\end{cases}$$
Since $J(0)=0$ and $J(tu)\longrightarrow -\infty $ as $t\longrightarrow +\infty$, then there exist a unique $t_u>0$ such that $t_uu\in \mathcal{N}$. \qed

\vspace{0.4cm}
\bl\label{lemma} Let $0 <\alpha <N$,~$ p= \dfrac{N+\alpha}{N-2}$, when ~$N\geqslant 5$,~$2< q < 2^{*}$, or $N=4$, $3<q<4$, we have  $0<c<\dfrac{p-1}{2p}C_0^{-\frac{1}{p-1}}S^{\frac{p}{p-1}}$. \el
\vspace{0.4cm}
\noindent {\bf Proof. } For $\forall u\in \mathcal{N}$, by Hardy-Littlewood-Sobolev inequality and Sobolev embedding theorem , it's easy to verify that there exists a  $\delta_0>0$ such that 
$$\rVert u \lVert \geqslant \delta_0 .$$
Then we obtain that 
\begin{equation}
\begin{split}
I(u)&=\dfrac{1}{2}[A(u)+D(u)]-\dfrac{1}{2p}B(u)-\dfrac{1}{q}C(u)\\
&\geqslant \big(\dfrac{1}{2}-\max\{\dfrac{1}{2p},\dfrac{1}{q}\}\big)\rVert u \lVert^2\\
&\geqslant \big(\dfrac{1}{2}-\max\{\dfrac{1}{2p},\dfrac{1}{q}\}\big)\delta_0^2\\
&>0.
\end{split}
\end{equation}
From the definition of $c$, we get that $c>0$. 

To prove $c<\dfrac{p-1}{2p}C_0^{-\frac{1}{p-1}}S^{\frac{p}{p-1}}$, firstly we consider the case of $N\geqslant 5$. Set
$$U(x)=\dfrac{[N(N-2)]^{(N-2)/4}}{[1+|x|^2]^{(N-2)/2}},$$
$$U_\varepsilon(x)=\varepsilon^\frac{2-N}{2}U(x/\varepsilon)=\varepsilon^\frac{2-N}{2}\dfrac{[N(N-2)]^{(N-2)/4}}{[1+|x|^2/\varepsilon^2]^{(N-2)/2}},$$
as in [27], then we have $$\lVert \nabla U_\varepsilon \rVert_2^2= \lVert U_\varepsilon \rVert_{2^*}^{2^*}=S^{N/2}.$$
By Prop 2.1, we have that 
$$B(U_\varepsilon)=C_0\lVert U_\varepsilon \rVert_{2^*}^{2p}=C_0S^{(N+\alpha)/2}.$$ 
Also, 
\begin{equation}
\begin{split}
\lVert  U_\varepsilon \rVert_2^2&=\varepsilon^{2-N}\int_{\mathbb{R}^N}\dfrac{[N(N-2)]^{(N-2)/2}}{[1+|x|^2/\varepsilon^2]^{N-2}}dx\\
&=C\varepsilon^2\int_{\mathbb{R}^N}\dfrac{1}{[1+|y|^2]^{N-2}}dy\\
&=C\varepsilon^2\int_{0}^{+\infty} \dfrac{r^{N-1}}{[1+r^2]^{N-2}}dr\\
&=O(\varepsilon^2),
\end{split}
\end{equation}
and
\begin{equation}
\begin{split}
\lVert  U_\varepsilon \rVert_q^q&=\varepsilon^{q(2-N)/2}\int_{\mathbb{R}^N}\dfrac{[N(N-2)]^{q(N-2)/4}}{[1+|x|^2/\varepsilon^2]^{q(N-2)/2}}dx\\
&=C\varepsilon^\frac{{2N-q(N-2)}}{2}\int_{\mathbb{R}^N}\dfrac{1}{[1+|y|^2]^{q(N-2)/2}}dy\\
&=C\varepsilon^\frac{{2N-q(N-2)}}{2}\int_{0}^{+\infty} \dfrac{r^{N-1}}{[1+r^2]^{q(N-2)/2}}dr\\
&=O(\varepsilon^\frac{{2N-q(N-2)}}{2}).
\end{split}
\end{equation}
Since 
\begin{equation}
\begin{split}
c=\inf\limits_{u\in \mathcal{N}}I(u)&\leqslant \max\limits_{t\geqslant 0}I(t U_\varepsilon)\\
&=\max\limits_{t\geqslant 0} [\dfrac{t^2}{2}A( U_\varepsilon)+\dfrac{t^2}{2}D( U_\varepsilon)-\dfrac{t^{2p}}{2p}B( U_\varepsilon)-\dfrac{t^q}{q}C( U_\varepsilon)]\\
&=\dfrac{t_\varepsilon^2}{2}A( U_\varepsilon)+\dfrac{t_\varepsilon^2}{2}D( U_\varepsilon)-\dfrac{t_\varepsilon^{2p}}{2p}B( U_\varepsilon)-\dfrac{t_\varepsilon^q}{q}C( U_\varepsilon)\\
&\leqslant \max\limits_{t\geqslant 0} \Big(\dfrac{t^2}{2}A( U_\varepsilon)-\dfrac{t^{2p}}{2p}B( U_\varepsilon)\Big)+\dfrac{t_\varepsilon^2}{2}D( U_\varepsilon)-\dfrac{t_\varepsilon^q}{q}C( U_\varepsilon)\\
&=\dfrac{p-1}{2p}\big(\dfrac{A( U_\varepsilon)^p}{B( U_\varepsilon)}\big)^{\frac{1}{p-1}}+\dfrac{t_\varepsilon^2}{2}D( U_\varepsilon)-\dfrac{t_\varepsilon^q}{q}C( U_\varepsilon)\\
&=\dfrac{p-1}{2p}C_0^{-\frac{1}{p-1}}S^{\frac{p}{p-1}}+\dfrac{t_\varepsilon^2}{2}D( U_\varepsilon)-\dfrac{t_\varepsilon^q}{q}C( U_\varepsilon),
\end{split}
\end{equation}
choosing $\varepsilon >0$ small enough, we know that $t_\varepsilon$ is close to $C_0^{-\frac{1}{2(p-1)}}S^{-\frac{\alpha}{4(p-1)}}$, so we have that $c<\dfrac{p-1}{2p}C_0^{-\frac{1}{p-1}}S^{\frac{p}{p-1}}$. 

Now we consider the case of $N=4$. For $\varepsilon,\sigma >0$ , set
$$U^\sigma(x)=\dfrac{[N(N-2)]^{(N-2)/4}}{[1+|x|^2]^{(N-2+\sigma)/2}},$$
$$U_\varepsilon^\sigma(x)=\varepsilon^\frac{2-N}{2}U^\sigma(x/\varepsilon)=\varepsilon^\frac{2-N}{2}\dfrac{[N(N-2)]^{(N-2)/4}}{[1+|x|^2/\varepsilon^2]^{(N-2+\sigma)/2}}.$$
Choosing $\sigma=\varepsilon^s$, where $4-q<s<q-2$, for $\varepsilon >0$ small, we can get that 
\begin{equation}
\begin{split}
\lVert  U_\varepsilon^\sigma \rVert_2^2&=\varepsilon^{2-N}\int_{\mathbb{R}^N}\dfrac{[N(N-2)]^{(N-2)/2}}{[1+|x|^2/\varepsilon^2]^{N-2+\sigma}}dx\\
&=C\varepsilon^2\int_{\mathbb{R}^N}\dfrac{1}{[1+|y|^2]^{N-2+\sigma}}dy\\
&=C\varepsilon^2\int_{0}^{+\infty} \dfrac{r^{N-1}}{[1+r^2]^{N-2+\sigma}}dr\\
&=O(\varepsilon^{2-s}),
\end{split}
\end{equation}
and
\begin{equation}
\begin{split}
\lVert  U_\varepsilon^\sigma \rVert_q^q&=\varepsilon^{q(2-N)/2}\int_{\mathbb{R}^N}\dfrac{[N(N-2)]^{q(N-2)/4}}{[1+|x|^2/\varepsilon^2]^{q(N-2+\sigma)/2}}dx\\
&=C\varepsilon^{4-q}\int_{\mathbb{R}^N}\dfrac{1}{[1+|y|^2]^{q(N-2+\sigma)/2}}dy\\
&=C\varepsilon^{4-q}\int_{0}^{+\infty} \dfrac{r^{N-1}}{[1+r^2]^{q(N-2+\sigma)/2}}dr\\
&=O(\varepsilon^{4-q}).
\end{split}
\end{equation}
It is easy to verify that 
$$A(U_\varepsilon^\sigma)=A(U^\sigma),~~B(U_\varepsilon^\sigma)=B(U^\sigma).$$
With simple calculations,  for $\varepsilon >0$ small, we get that  
\begin{equation}
\begin{split}
A(U^\sigma)&=[N(N-2)]^{(N-2)/2}\int_{\mathbb{R}^N}\dfrac{(N-2+\sigma)^2|x|^2}{(1+|x|^2)^{N+\sigma}}dx\\
&\leqslant [N(N-2)]^{(N-2)/2}\int_{\mathbb{R}^N}\dfrac{(N-2)^2|x|^2}{(1+|x|^2)^N}dx+C\sigma\\
&=A(U)+C\sigma,
\end{split}
\end{equation}
and by Hardy-Littlewood-Sobolev inequality,
\begin{equation}
\begin{split}
B(U)-B(U^\sigma)=&C\int_{\mathbb{R}^N}\int_{\mathbb{R}^N}\dfrac{1}{|x-y|^{N-\alpha}(1+|x|^2)^p(1+|y|^2)^p}\\
&\cdot \big\{ 1-\dfrac{1}{[(1+|x|^2)(1+|y|^2)]^{p\sigma/2}}\big\}dxdy\\
\leqslant &C\int_{\mathbb{R}^N}\int_{\mathbb{R}^N}\dfrac{\sigma \ln[(1+|x|^2)(1+|y|^2)]}{|x-y|^{N-\alpha}(1+|x|^2)^p(1+|y|^2)^p}dxdy\\
\leqslant& C\sigma\int_{\mathbb{R}^N}\int_{\mathbb{R}^N}\dfrac{1}{|x-y|^{N-\alpha}(1+|x|^2)^{p-1}(1+|y|^2)^{p-1}}dxdy\\
\leqslant& C\sigma\big(\int_{\mathbb{R}^N}(\dfrac{1}{1+|x|^2})^{\frac{2N(p-1)}{N+\alpha}}dx\big)^\frac{N+\alpha}{N}\\
=&C\sigma,
\end{split}
\end{equation}
then we have that $A(U^\sigma)^p\leqslant A(U)^p+C\sigma$, $B(U^\sigma)\geqslant B(U)-C\sigma$, so $\dfrac{A( U^\sigma)^p}{B( U^\sigma)}\leqslant \dfrac{A( U)^p}{B( U)}+C\sigma$. Thus, $\big(\dfrac{A( U^\sigma)^p}{B( U^\sigma)}\big)^{\frac{1}{p-1}}\leqslant \big(\dfrac{A( U)^p}{B( U)}\big)^{\frac{1}{p-1}}+C\sigma$.
Then we obtain that
\begin{equation}
\begin{split}
c=\inf\limits_{u\in \mathcal{N}}I(u)&\leqslant \max\limits_{t\geqslant 0}I(t U_\varepsilon^\sigma)\\
&=\max\limits_{t\geqslant 0} [\dfrac{t^2}{2}A( U_\varepsilon^\sigma)+\dfrac{t^2}{2}D( U_\varepsilon^\sigma)-\dfrac{t^{2p}}{2p}B( U_\varepsilon^\sigma)-\dfrac{t^q}{q}C( U_\varepsilon^\sigma)]\\
&=\dfrac{t_\varepsilon^2}{2}A( U_\varepsilon^\sigma)+\dfrac{t_\varepsilon^2}{2}D( U_\varepsilon^\sigma)-\dfrac{t_\varepsilon^{2p}}{2p}B( U_\varepsilon^\sigma)-\dfrac{t_\varepsilon^q}{q}C( U_\varepsilon^\sigma)\\
&\leqslant \max\limits_{t\geqslant 0} [\dfrac{t^2}{2}A( U^\sigma)-\dfrac{t^{2p}}{2p}B( U^\sigma)]+\dfrac{t_\varepsilon^2}{2}D( U_\varepsilon^\sigma)-\dfrac{t_\varepsilon^q}{q}C( U_\varepsilon^\sigma)\\
&=\dfrac{p-1}{2p}\big(\dfrac{A( U^\sigma)^p}{B( U^\sigma)}\big)^{\frac{1}{p-1}}+\dfrac{t_\varepsilon^2}{2}D( U_\varepsilon^\sigma)-\dfrac{t_\varepsilon^q}{q}C( U_\varepsilon^\sigma)\\
&\leqslant\dfrac{p-1}{2p}C_0^{-\frac{1}{p-1}}S^{\frac{p}{p-1}}+O(\varepsilon^s)+\dfrac{t_\varepsilon^2}{2}O(\varepsilon^{2-s})-\dfrac{t_\varepsilon^q}{q}O(\varepsilon^{4-q})
\end{split}
\end{equation}
for $\varepsilon >0$ small enough, and we have that $c<\dfrac{p-1}{2p}C_0^{-\frac{1}{p-1}}S^{\frac{p}{p-1}}$.
\qed
\vspace{0.4cm}
\bl\label{lemma} The $(PS)_c$ sequence of the constrained functional $I|_\mathcal{N}$ is also a $(PS)_c$ sequence of $I$ , namely,  if $\{u_n \}\subset \mathcal{N}$ satisfies $I(u_n)\longrightarrow c$ and $I'|_\mathcal{N}(u_n)\longrightarrow 0$, then $I'(u_n)\longrightarrow 0$. \el
\vspace{0.4cm}
\noindent {\bf Proof. } Firstly we prove $\{u_n \}$ is bounded. Since $\{u_n \}\subset \mathcal{N}$, then $A(u_n)+D(u_n)=B(u_n)+C(u_n)$. Since $I(u_n)\longrightarrow c$, we have 
\begin{equation}
\begin{split}
c+o(1)=I(u_n)&=\dfrac{1}{2}\lVert u_n \rVert^2-\dfrac{1}{2p}B(u_n)-\dfrac{1}{q}C(u_n)\\
&\geqslant \dfrac{1}{2}\lVert u_n \rVert^2-\max\{\dfrac{1}{2p},\dfrac{1}{q}\}[B(u_n)+C(u_n)]\\
&=\big(\dfrac{1}{2}-\max\{\dfrac{1}{2p},\dfrac{1}{q}\}\big)\lVert u_n \rVert^2,
\end{split}
\end{equation}
thus $\{u_n \}$ is bounded.

By Lagrange multiplier theorem, there exists a sequence $\{t_n\}\subset \mathbb{R}$ such that 
$$I'(u_n)=I'|_\mathcal{N}(u_n)+t_nJ'(u_n).$$
Then we have $$0=\left<I'(u_n),u_n\right>=\left<I'|_\mathcal{N}(u_n),u_n\right>+t_n\left<J'(u_n),u_n\right>,$$
therefore $$t_n\left<J'(u_n),u_n\right>\longrightarrow 0.$$
Since 
\begin{equation}
\begin{split}
\left<J'(u_n),u_n\right>&=2A(u_n)+2D(u_n)-2pB(u_n)-qC(u_n)\\
&\leqslant 2\lVert u_n \rVert^2 -\min\{2p,q\}\big(B(u_n)+C(u_n)\big)\\
&\leqslant\big(2-\min\{2p,q\}\big)\delta_0^2\\
&<0,
\end{split}
\end{equation}
and 
$$\left<J'(u_n),u_n\right>\geqslant \big(2-\max\{2p,q\}\big)\lVert u_n \rVert^2,$$
then we get that $\{\left<J'(u_n),u_n\right>\}$ is bounded due to the boundedness of $\{u_n\}$. So $t_n\longrightarrow 0$, and $I'(u_n)\longrightarrow 0$. \qed

\vspace{0.4cm}
\noindent {\bf Proof of Theorem1.1. } By Ekeland variational principle, there exists a sequence $\{u_n \}\subset \mathcal{N}$ satisfying $I(u_n)\longrightarrow c$ and $I'|_\mathcal{N}(u_n)\longrightarrow 0$. By Lemma 2.3, $\{u_n \}$ is bounded and $I'(u_n)\longrightarrow 0$. Then for $\forall v \in  H^1(\mathbb{R}^N)$, we have $$\int_{\mathbb{R}^N}\nabla u_n\cdot \nabla v+\int_{\mathbb{R}^N}u_nv-\int_{\mathbb{R}^N}(I_\alpha*|u_n|^p)|u_n|^{p-2}u_nv-\int_{\mathbb{R}^N}|u_n|^{q-2}u_nv=o(1)\lVert v \rVert.$$
Assume $u_n \rightharpoonup u_0$ up to subsequence, and $u_n \longrightarrow u_0$ a.e. ~on $\mathbb{R}^N$ . Then we have 
$$\int_{\mathbb{R}^N}\nabla u_n\cdot \nabla v \longrightarrow \int_{\mathbb{R}^N}\nabla u_0\cdot \nabla v,$$
$$\int_{\mathbb{R}^N}u_nv\longrightarrow \int_{\mathbb{R}^N}u_0v,$$
$$\int_{\mathbb{R}^N}|u_n|^{q-2}u_nv\longrightarrow \int_{\mathbb{R}^N}|u_0|^{q-2}u_0v.$$
Moreover, for $\forall \phi \in C_c^1( \mathbb{R}^N)$, suppose that $\mbox{supp}(\phi)\subset K\subset \subset \mathbb{R}^N$, by Rellich compactness theorem and Sobolev inequality, we get that $|u_n|^p\rightharpoonup |u_c|^p$ in $L^{\frac{2N}{N+\alpha}}(\mathbb{R}^N)$ and $|u_n|^{p-2}u_n \longrightarrow |u_c|^{p-2}u_c$ in measure on $K$ up to subsequence. By Hardy-Littlewood-Sobolev inequality, we have that $I_\alpha*|u_n|^p\rightharpoonup I_\alpha*|u_c|^p$ in $L^{\frac{2N}{N-\alpha}}(\mathbb{R}^N)$. Since $|u_n|^{p-2}u_n \longrightarrow |u_c|^{p-2}u_c$ in measure on $K$, for $\forall \varepsilon>0, \delta>0$, $\exists N>0$  such that $\forall n>N$, we have 
$$m\big\{x\in K:\big| |u_n(x)|^{p-2}u_n(x)-|u_c(x)|^{p-2}u_c(x)\big|>\delta\big\}<\varepsilon.$$
Set $K_1=\big\{x\in K:\big| |u_n(x)|^{p-2}u_n(x)-|u_c(x)|^{p-2}u_c(x)\big|>\delta\big\}$. By Hardy-Littlewood-Sobolev inequality and H\"{o}lder inequality, we have  

\begin{equation}   \nonumber
\begin{split}
&\big|\int_{\mathbb{R}^N}(I_\alpha*|u_n|^p)|u_n|^{p-2}u_n\phi- \int_{\mathbb{R}^N}(I_\alpha*|u_c|^p)|u_c|^{p-2}u_c\phi\big|\\
&\leqslant \Big|\int_{\mathbb{R}^N}\big(I_\alpha*(|u_n|^p-|u_c|^p)\big)|u_c|^{p-2}u_c\phi\Big|+\int_{K}(I_\alpha*|u_n|^p)\big||u_n|^{p-2}u_n-|u_c|^{p-2}u_c\big||\phi|\\
&:=I_1+I_2.
\end{split}
\end{equation}
Since $I_\alpha*|u_n|^p\rightharpoonup I_\alpha*|u_c|^p$ in $L^{\frac{2N}{N-\alpha}}(\mathbb{R}^N)$, then $I_1\longrightarrow 0$. And
\begin{equation}   \nonumber
\begin{split}
I_2=&\int_{K_1}(I_\alpha*|u_n|^p)\big||u_n|^{p-2}u_n-|u_c|^{p-2}u_c\big||\phi|+\int_{K/K_1}(I_\alpha*|u_n|^p)\big||u_n|^{p-2}u_n-|u_c|^{p-2}u_c\big||\phi|\\
\leqslant& C\big(\int_{\mathbb{R}^N}|u_n|^\frac{2Np}{N+\alpha}\big)^{\frac{N+\alpha}{2N}}\Big(\big(\int_{K_1}\big||u_n|^{p-2}u_n-|u_c|^{p-2}u_c\big|^\frac{2Np}{(N+\alpha)(p-1)}\big)^{\frac{(N+\alpha)(p-1)}{2Np}}\big(\int_{K_1}|\phi|^\frac{2Np}{N+\alpha}\big)^\frac{N+\alpha}{2Np}\\
&+\delta\big(\int_{K/K_1}|\phi|^\frac{2Np}{N+\alpha}\big)^\frac{N+\alpha}{2Np}\Big),
\end{split}
\end{equation}
by the absolute continuity of the integral,  we have $$\int_{\mathbb{R}^N}(I_\alpha*|u_n|^p)|u_n|^{p-2}u_n\phi\longrightarrow \int_{\mathbb{R}^N}(I_\alpha*|u_c|^p)|u_c|^{p-2}u_c\phi.$$
Since $C^1_c(\mathbb{R}^N)$ is dense in $H^1(\mathbb{R}^N)$,then we have 
$$  \int_{\mathbb{R}^N}(I_\alpha*|u_n|^p)|u_n|^{p-2}u_nv\longrightarrow \int_{\mathbb{R}^N}(I_\alpha*|u_0|^p)|u_0|^{p-2}u_0v.$$
Thus we get that $I'(u_0)=0$.

Finally, we want to show that $u_0\neq 0$. Otherwise, we assume that $u_n\rightharpoonup u_0=0$. Since the embedding $H_{rad}^1(\mathbb{R}^N)\hookrightarrow L^q(\mathbb{R}^N)$ is compact,
 we obtain that $C(u_n)\longrightarrow 0$. Then we have 
$$\lVert u_n \rVert^2\longrightarrow \dfrac{2p}{p-1}c,$$
$$B(u_n)\longrightarrow \dfrac{2p}{p-1}c.$$
By Hardy-Littlewood-Sobolev inequality, we know that 
$$B(u_n)\leqslant C_0S^{-p}\lVert \nabla u_n\rVert_2^{2p}\leqslant C_0S^{-p}\lVert  u_n\rVert^{2p}.$$
Letting $n \longrightarrow +\infty$, we have that 
$$\dfrac{2p}{p-1}c\leqslant C_0S^{-p}(\dfrac{2p}{p-1}c)^p,$$
then 
$$c\geqslant \dfrac{p-1}{2p}C_0^{-\frac{1}{p-1}}S^{\frac{p}{p-1}},$$
which contradicts with Lemma 2.2. So we get that $u_0\neq 0$, and $u_0$ is a solution  of (1.1).\qed

\vspace{0.4cm}
\section{The Proof of Theorem 1.2}

\vspace{0.4cm}
Define the set
$$\mathcal{M}:=\{u \in H^1(\mathbb{R}^N)\setminus \{0\}:\int_{\mathbb{R}^N}G(u)dx=1\}.$$
It is important to observe that $\mathcal{M}$ is a $C^1$ manifold. Indeed, for every $u\in \mathcal{M} $, let $H(u)=\int_{\mathbb{R}^N}G(u)dx$, then 
$$\left<H'(u),u\right>=B(u)+C(u)-D(u)=(1-\dfrac{1}{p})B(u)+(1-\dfrac{2}{q})C(u)+2>2,$$
thus $H'(u)\neq 0$.
\vspace{0.4cm}
\bl\label{lemma} Any minimizing sequence $\{u_n\}$ by (1.5) is bounded in $H^1(\mathbb{R}^N)$.\el
\vspace{0.4cm}
\noindent {\bf Proof. } Since $\{u_n\}$ is a minimizing sequence, we have
$$\frac{1}{2}A(u_n)\longrightarrow A,$$
and 
$$\int_{\mathbb{R}^N}G(u_n)dx=\frac{1}{2p}B(u_n)+\dfrac{1}{q}C(u_n)-\dfrac{1}{2}D(u_n)=1.$$
By Hardy-Littlewood-Sobolev inequality, we have
$$B(u_n)\leqslant C_0\parallel u_n\parallel_{2^*}^{2p},$$
where $C_0$ is the best constant. Moreover, there exists a constant $C>0$ such that 
$$C(u_n)=\int_{\mathbb{R}^N}|u_n|^qdx\leqslant \int_{\mathbb{R}^N}(\frac{1}{4}|u_n|^2+C|u_n|^{2^*})dx.$$
Then we have 
\begin{equation}
\begin{split}
\frac{1}{2}D(u_n)+1&=\frac{1}{2p}B(u_n)+\dfrac{1}{q}C(u_n){}\\
&\leqslant C(\lVert u_n\rVert_{2^*}^{2p}+\lVert u_n\rVert_{2^*}^{2^*})+\frac{1}{4}D(u_n).
\end{split}
\end{equation}
Thus, 
$$D(u_n)\leqslant  C(\lVert u_n\rVert_{2^*}^{2p}+\lVert u_n\rVert_{2^*}^{2^*})\leqslant C \big(S^{-\frac{N+\alpha}{N-2}}A(u_n)^{\frac{N+\alpha}{N-2}}+S^{-\frac{N}{N-2}}A(u_n)^{\frac{N}{N-2}}
\big).$$
Since $A(u_n)\longrightarrow 2A$, then $\{A(u_n)\}$ and $\{D(u_n)\}$ are both bounded, which implies the boundedness of $\{u_n\}$ in $H^1(\mathbb{R}^N)$.  \qed

\vspace{0.4cm}
\bl\label{lemma} Under the assumptions of theorem 1.2, we have that $0<A<\frac{1}{2}S\left(\dfrac{2p}{C_0}\right)^{\frac{1}{p}}$.\el
\vspace{0.4cm}
\noindent {\bf Proof. }  Firstly we want to prove the set $\mathcal{M}$ is not empty. When $N\geqslant 5$, set
$$u_\varepsilon(x)=\dfrac{\left(N(N-2)\varepsilon^2\right)^{\frac{N-2}{4}}}{(\varepsilon^2+|x|^2)^{\frac{N-2}{2}}}.$$ 
As in [27], we obtain
$$\lVert\nabla u_\varepsilon\rVert_2^2=\lVert  u_\varepsilon\rVert_{2^*}^{2^*}=S^{\frac{N}{2}}.$$
Moreover, we know that $u_\varepsilon$ is a maximizing function of the Hardy-Littlewood-Sobolev inequality, see [14] .
Let $v_\varepsilon=\dfrac{u_\varepsilon}{\lVert  u_\varepsilon\rVert_{2^*}}$, we have
\begin{equation}
T(v_\varepsilon)=\dfrac{\lVert\nabla u_\varepsilon\rVert_2^2}{2\lVert  u_\varepsilon\rVert_{2^*}^2}=\dfrac{S}{2},
\end{equation}

\begin{equation}
B(v_\varepsilon)=\dfrac{B(u_\varepsilon)}{\lVert  u_\varepsilon\rVert_{2^*}^{2p}}= \dfrac{C_0\lVert  u_\varepsilon\rVert_{2^*}^{2p}}{\lVert  u_\varepsilon\rVert_{2^*}^{2p}}=C_0,
\end{equation}

\begin{equation}
C(v_\varepsilon)=\dfrac{\int_{\mathbb{R}^N}|u_\varepsilon|^qdx}{\lVert  u_\varepsilon\rVert_{2^*}^q}=C_1\varepsilon^{N-\frac{N-2}{2}q}\int_{0}^{\infty}\dfrac{r^{N-1}}{(1+r^2)^{\frac{N-2}{2}q}}dr,
\end{equation}

\begin{equation}
D(v_\varepsilon)=\dfrac{\int_{\mathbb{R}^N}|u_\varepsilon|^2dx}{\lVert  u_\varepsilon\rVert_{2^*}^2}=C_2\varepsilon^2\int_{0}^{\infty}\dfrac{r^{N-1}}{(1+r^2)^{N-2}}dr.
\end{equation}

It is easy to check that $\int_{0}^{\infty}\dfrac{r^{N-1}}{(1+r^2)^{\frac{N-2}{2}q}}dr<+\infty$ and $\int_{0}^{\infty}\dfrac{r^{N-1}}{(1+r^2)^{N-2}}dr<+\infty$. 
For $t>0$, we have that 
$$H(tv_\varepsilon)=\dfrac{C_0}{2p}t^{2p}+C_1t^q\varepsilon^{N-\frac{N-2}{2}q}-C_2t^2\varepsilon^2.$$
Set $$g_\varepsilon(t)=C_1t^q\varepsilon^{N-\frac{N-2}{2}q}-C_2t^2\varepsilon^2.$$
Then there exists $\varepsilon_0>0$ small, such that for $\forall t\in \big[(\dfrac{p}{C_0})^\frac{1}{2p},(\dfrac{2p}{C_0})^\frac{1}{2p}\big]$, when $\varepsilon \in (0,\varepsilon_0)$, we have that 
$$0<g_\varepsilon(t)<\dfrac{1}{2}.$$
Then $\exists t_\varepsilon \in \big[(\dfrac{p}{C_0})^\frac{1}{2p},(\dfrac{2p}{C_0})^\frac{1}{2p}\big]$ such that $H(t_\varepsilon v_\varepsilon)=1$ because $H(tv_\varepsilon)$ is continuous in $t$. Thus the set $\mathcal{M}$ is not empty. 
Moreover, since
$$1=H(t_\varepsilon v_\varepsilon)=\dfrac{C_0}{2p}t_\varepsilon^{2p}+C_1t_\varepsilon^q\varepsilon^{N-\frac{N-2}{2}q}-C_2t_\varepsilon^2\varepsilon^2,$$
then $\dfrac{C_0}{2p}t_\varepsilon^{2p}<1$ for $\varepsilon \in (0,\varepsilon_0)$, from which we know that 
$$A\leqslant T(t_\varepsilon v_\varepsilon)=t_\varepsilon^2T(v_\varepsilon)=\dfrac{1}{2}St_\varepsilon^2<\dfrac{1}{2}S(\dfrac{2p}{C_0})^{\frac{1}{p}}.$$
When $N=4$, as in Lemma 2.2, set
$$u^\sigma(x)=\dfrac{[N(N-2)]^{(N-2)/4}}{[1+|x|^2]^{(N-2+\sigma)/2}},$$
$$u_\varepsilon^\sigma(x)=\varepsilon^\frac{2-N}{2}u^\sigma(x/\varepsilon)=\varepsilon^\frac{2-N}{2}\dfrac{[N(N-2)]^{(N-2)/4}}{[1+|x|^2/\varepsilon^2]^{(N-2+\sigma)/2}}.$$ 
Choosing $\sigma=\varepsilon^s$, where $4-q<s<q-2$, for $\varepsilon >0$ small, we have that
$$\lVert  u_\varepsilon^\sigma \rVert_2^2=O(\varepsilon^{2-s}),~~
\lVert  u_\varepsilon^\sigma \rVert_q^q=O(\varepsilon^{4-q}),$$
and $$A(u_\varepsilon^\sigma)\leqslant S^{N/2}+C\sigma,$$ $$C_0S^{(N+\alpha)/2}-C\sigma \leqslant B(u_\varepsilon^\sigma)\leqslant C_0S^{(N+\alpha)/2},$$ 
$$S^{N/2}-C\sigma\leqslant \lVert  u_\varepsilon^\sigma \rVert_{2^*}^{2^*}\leqslant S^{N/2} .$$ 
 Set $v_\varepsilon^\sigma=\dfrac{u_\varepsilon^\sigma}{\lVert  u_\varepsilon^\sigma\rVert_{2^*}}$,  we have that 
\begin{equation}
T(v_\varepsilon^\sigma)=\dfrac{\lVert\nabla u_\varepsilon^\sigma\rVert_2^2}{2\lVert  u_\varepsilon^\sigma\rVert_{2^*}^2}\leqslant \dfrac{S}{2}+C\sigma,
\end{equation}

\begin{equation}
C_0-C\sigma \leqslant B(v_\varepsilon^\sigma)=\dfrac{B(u_\varepsilon^\sigma)}{\lVert  u_\varepsilon^\sigma\rVert_{2^*}^{2p}}\leqslant C_0+C\sigma,
\end{equation}

\begin{equation}
O(\varepsilon^{4-q}) \leqslant C(v_\varepsilon^\sigma)=\dfrac{\int_{\mathbb{R}^N}|u_\varepsilon^\sigma|^qdx}{\lVert  u_\varepsilon^\sigma\rVert_{2^*}^q}\leqslant O(\varepsilon^{4-q})+C\sigma,
\end{equation}

\begin{equation}
O(\varepsilon^{2-s}) \leqslant D(v_\varepsilon^\sigma)=\dfrac{\int_{\mathbb{R}^N}|u_\varepsilon^\sigma|^2dx}{\lVert  u_\varepsilon^\sigma\rVert_{2^*}^2}\leqslant O(\varepsilon^{2-s})+C\sigma.
\end{equation}
For $\forall t\in \big[(\dfrac{p}{C_0})^\frac{1}{2p},(\dfrac{2p}{C_0})^\frac{1}{2p}\big]$, we have that 
$$\Big|H(tv_\varepsilon^\sigma)-\big(\dfrac{C_0}{2p}t^{2p}+\dfrac{t^q}{q}O(\varepsilon^{4-q})-\dfrac{t^2}{2}O(\varepsilon^{2-s})\big)\Big|<C\sigma,$$
and $T(v_\varepsilon^\sigma)\leqslant \dfrac{S}{2}+C\sigma<\dfrac{S}{2}(1+C_1\sigma)^{\frac{1}{p}}$.
Set $$g_\varepsilon^\pm(t)=\dfrac{t^q}{q}O(\varepsilon^{4-q})-\dfrac{t^2}{2}O(\varepsilon^{2-s})\pm C\sigma.$$
Then there exists $\varepsilon_0>0$ small, such that for $\forall t\in \big[(\dfrac{p}{C_0})^\frac{1}{2p},(\dfrac{2p}{C_0})^\frac{1}{2p}\big]$, when $\varepsilon \in (0,\varepsilon_0)$, we have that 
$$0<g_\varepsilon^\pm(t)<\dfrac{1}{2}.$$
Then $\exists t_\varepsilon \in \big[(\dfrac{p}{C_0})^\frac{1}{2p},(\dfrac{2p}{C_0})^\frac{1}{2p}\big]$ such that $H(t_\varepsilon v_\varepsilon^\sigma)=1$ because $H(tv_\varepsilon^\sigma)$ is continuous in $t$. Thus the set $\mathcal{M}$ is not empty. Moreover, since
$$1=H(t_\varepsilon v_\varepsilon^\sigma)>\dfrac{C_0}{2p}t_\varepsilon^{2p}+t_\varepsilon^qO(\varepsilon^{4-q})-t_\varepsilon^2O(\varepsilon^{2-s})-C\sigma,$$
then $\dfrac{C_0}{2p}t_\varepsilon^{2p}<1-C_1\sigma$ for $\varepsilon>0$ small enough, from which we know that 
$$A\leqslant T(t_\varepsilon v_\varepsilon^\sigma)=t_\varepsilon^2T(v_\varepsilon^\sigma)<\dfrac{1}{2}S(1+C_1\sigma)^{\frac{1}{p}}(1-C_1\sigma)^{\frac{1}{p}}(\dfrac{2p}{C_0})^{\frac{1}{p}}<\dfrac{1}{2}S(\dfrac{2p}{C_0})^{\frac{1}{p}}.$$
On the other hand, as in the proof of Lemma3.1, we have that $A>0$.
\qed

\vspace{0.4cm}
The next lemma is the Brezis-Lieb lemma for the nonlocal term of the functional, see [21].
\vspace{0.4cm}
\bl\label{lemma} Let $N\in \mathbb{N},~\alpha\in (0,N),~p\in [1,\frac{N+\alpha}{N-2}]$ and $\{u_n\}$ be a bounded sequence in $L^{2Np/(N+\alpha)}(\mathbb{R}^N)$. If $u_n\longrightarrow u $ a.e. on $\mathbb{R}^N$, then 
$$\int_{\mathbb{R}^N}(I_\alpha*|u_n|^p)|u_n|^pdx-\int_{\mathbb{R}^N}(I_\alpha*|u_n-u|^p)|u_n-u|^pdx=\int_{\mathbb{R}^N}(I_\alpha*|u|^p)|u|^pdx+o(1).$$\el

\vspace{0.4cm}
\bl\label{lemma} Assume the assumptions of Theorem1.2 hold. Then problem (1.5)
has a minimizer $u_0\in H_{rad}^1(\mathbb{R}^N)$. \el
\vspace{0.4cm}
\noindent {\bf Proof. } First we show there exists a radial minimizing sequence. Assume $\{u_n\}\in H^1(\mathbb{R}^N)$ is a minimizing sequence , namely $H(u_n)=1$, and $T(u_n)\longrightarrow A$. Let $u_n^*$ be the Schwarz spherical rearrangement of $|u_n|$, then $u_n^*$ is radial and by [15], we have that $H(u_n^*)\geqslant H(u_n)=1$, $T(u_n^*)\leqslant T(u_n)$. Then there exists $t_n\in [0,1]$, such that $H(t_nu_n^*)=1$. Set $\bar{u}_n=t_nu_n^*$, then $H(\bar{u}_n)=1$, $T(\bar{u}_n)=t_n^2T(u_n^*)\leqslant T(u_n^*)$, so $\{\bar{u}_n\}$ is a  radial minimizing sequence. So we assume $\{u_n\}\in H_{rad}^1(\mathbb{R}^N)$ is a minimizing sequence. By Lemma 2.1, we know that  $\{u_n\}$ is bounded in $H_{rad}^1(\mathbb{R}^N)$. We can assume $u_n\rightharpoonup u_0$ in $H^1(\mathbb{R}^N)$, and $u_n\longrightarrow u_0 ~a.e.~$ on $\mathbb{R}^N$ up to subsequence. Then we have $T(u_0)\leqslant A$. By Radial Lemma in [26], we know that $|u_n(x)|\leqslant C|x|^{-\frac{N}{2}},~ x\in \mathbb{R}^N$. Set $v_n=u_n-u_0$, then we have $T(u_n)=T(v_n)+T(u_0)+o(1)$. By Brezis-Lieb Lemma in [27], we get that 
$C(u_n)=C(v_n)+C(u_0)+o(1)$ and $D(u_n)=D(v_n)+D(u_0)+o(1)$. From Lemma 3.3, we know that $B(u_n)=B(v_n)+B(u_0)+o(1)$. Thus we obtain $H(u_n)=H(v_n)+H(u_0)+o(1)$. Set $S_n=T(v_n)$, $S_0=T(u_0)$, $\lambda_n=H(v_n)$, $\lambda_0=H(u_0)$, then we have that 
$$\lambda_n=1-\lambda_0+o(1),~~S_n=A-S_0+o(1).$$
It is sufficient to prove that $\lambda_0=1$. Set $u_\sigma(x)=u(\dfrac{x}{\sigma}),~\sigma>0$, then we have 
$$T(u_\sigma)=\sigma^{N-2}T(u),~B(u_\sigma)=\sigma^{N+\alpha}B(u),$$
$$C(u_\sigma)=\sigma^{N}C(u),~D(u_\sigma)=\sigma^{N}D(u).$$
Since $H(u_\sigma)=\frac{\sigma^{N+\alpha}}{2p}B(u)+\frac{\sigma^N}{q}C(u)-\frac{\sigma^N}{2}D(u)$, for $\forall u\neq 0$, there exists $\sigma_u>0$ such that $H(u_\sigma)=1$.
By the definition of $A$, we have that 
$$T(u_{\sigma_u})\geqslant A H(u_{\sigma_u})=A\big[{\sigma_u}^{N}H(u)+\dfrac{{\sigma_u}^N}{2p}({\sigma_u}^\alpha-1)B(u)\big].$$
When $0<H(u)\leqslant 1$, we have $\sigma_u\geqslant 1$, then ${\sigma_u}^{N}H(u)\leqslant 1$ and we know that $T(u)\geqslant A{\sigma_u}^{2-N}\geqslant AH(u)^{\frac{N-2}{N}}$. When $H(u) > 1$, we get that $\sigma_u<1$, then $T(u)\geqslant A{\sigma_u}^{2-N}>A$. If $\lambda_0>1$, then $S_0>A$, which is a contradiction with $S_0\leqslant A$. So $\lambda_0\leqslant 1$. If $\lambda_0<0$, then $\lambda_n>1-\dfrac{\lambda_0}{2}>1$ for $n$ large enough. Then there exists $\varepsilon_0>0$ such that $\sigma_{v_n}\leqslant 1-\varepsilon_0$ for $n$ large enough. Thus, we have $S_n\geqslant  A{\sigma_{v_n}}^{2-N}\geqslant A( 1-\varepsilon_0)^{2-N}$ for $n$ large enough , which is a contradiction to $S_n\leqslant A+o(1)$. Thus we get that $\lambda_0 \in [0,1]$. If $\lambda_0 \in (0,1)$, then $0<\lambda_n<1$ for n large enough, and $S_0\geqslant A(\lambda_0)^\frac{N-2}{N}$, $S_n\geqslant A(\lambda_n)^\frac{N-2}{N}$. Then we have that
\begin{equation}
\begin{split}
A&=\lim\limits_{n\longrightarrow \infty}(S_0+S_n)\\
&\geqslant \lim\limits_{n\longrightarrow \infty}A(\lambda_0^\frac{N-2}{N}+\lambda_n^\frac{N-2}{N})\\
&=A(\lambda_0^\frac{N-2}{N}+(1-\lambda_0)^\frac{N-2}{N})\\
&> A(\lambda_0+1-\lambda_0)=A,
\end{split}
\end{equation}
which is a contradiction.
If $\lambda_0=0$, then $\lambda_n=1+o(1)$, and we get that $S_0=0$ and $u_0=0$. Thus  
$$1=H(v_n)=\dfrac{1}{2p}B(v_n)+\dfrac{1}{q}C(v_n)-\dfrac{1}{2}D(v_n).$$
By Hardy-Littlewood-Sobolev inequality, we have that 
$$1 +\dfrac{1}{2}D(v_n)=\dfrac{1}{2p}B(v_n)+\dfrac{1}{q}C(v_n)\leqslant \dfrac{C_0}{2p}\lVert v_n\rVert_{2^*}^{2p}+\dfrac{1}{q}C(v_n).$$
Since the embedding $H_{rad}^1(\mathbb{R}^N)\subset L^q(\mathbb{R}^N)$ is compact, we have $C(v_n)=o(1)$. Then 
$$\limsup\limits_{n\longrightarrow \infty}\lVert v_n\rVert_{2^*}^2\geqslant (\dfrac{2p}{C_0})^{\frac{1}{p}},$$
then by Sobolev inequality,
\begin{equation}
\begin{split}
A&=\dfrac{1}{2}\lim_{n\longrightarrow \infty}\lVert \nabla v_n\rVert_2^2\\
&\geqslant \dfrac{1}{2}(\dfrac{2p}{C_0})^{\frac{1}{p}}\liminf\limits_{n\longrightarrow \infty}\dfrac{\lVert \nabla v_n\rVert_2^2}{\lVert v_n\rVert_{2^*}^2}\\
&\geqslant \dfrac{1}{2}(\dfrac{2p}{C_0})^{\frac{1}{p}}S,
\end{split}
\end{equation}
which is a contradiction by Lemma 3.2. Therefore, we conclude that $\lambda_0=1$, and $u_0$ is a minimizer. \qed

\vspace{0.4cm}
%\vskip0.7in
%\noindent{\bf References} \vskip0.2in
\footnotesize

%\bibliographystyle{abbrv}
%\bibliography{papers,books}

%%%%%%%%%%%%%%%%%%%%%%%%%%%%%%%%%%%%%%%%%%%%%%%%%%%%%%%%%%%%%%%%%%%%%%%%%%
\end{CJK*}
\end{document}